\documentclass{article}
\usepackage{amsmath,amsfonts,amsthm,amssymb,amscd}
 \usepackage{amssymb}
 \usepackage[cp1251]{inputenc}
 \usepackage[T2A]{fontenc}
 \usepackage[all]{xy}
 \usepackage{latexsym}
 \usepackage{srcltx}
 \usepackage{amsmath}
 \usepackage{amsfonts}
 \usepackage{amscd}
 \usepackage{amssymb}
 \usepackage{amsthm}
 \usepackage{euscript}
 \usepackage{url}
 \usepackage{amssymb, amsthm}
 \usepackage{mathrsfs}
 \usepackage{srcltx}
 \usepackage{hyperref}\usepackage{setspace}
 \usepackage[russian,ukrainian,english]{babel}
 \usepackage[pdftex]{graphicx}%
 \usepackage{mathrsfs}
 \usepackage{graphicx}

 \newtheorem{proposition}{Proposition}
 \newtheorem{theorem}{Theorem}

 \newtheorem{definition}{Definition}
\newtheorem{corollary}{Corollary}

\begin{document}


 \begin{center}
   \textbf{FIBRATIONS OF FINANCIAL EVENTS\\[.5cm]
       David Carf\`{i}}\\[.2cm]
   (Messina University, Messina, Italy)\\[.2cm]
   \emph{E-mail address}: \texttt{davidcarfi71@yahoo.it}\\[.3cm]
 \end{center}
 \bigskip

\begin{abstract}
In this paper we shall prove that the plane of financial events, introduced
and applied to financial problems by the author himself (see \cite{car2}, 
\cite{car3} and \cite{car5}) can be considered as a fibration in two
different ways. The first one, \emph{the natural one}, reveals itself to be
isomorphic to the tangent-bundle of the real line, when the last one is
considered as a differentiable manifold in the natural way; the second one
is a fibration induced by the \emph{status of compound interest
capitalization} at a given rate $i\in ]-1,\rightarrow [$. Moreover, in the
paper we define on the first fibration an affine connection, also in this
case induced by the status of compound interest at a given rate $i$. The
final goal of this paper is the awareness that all the effects determined by
the status of compound interest are \emph{nothing but} the consequences of
the fact that the space of financial events is a fibration endowed with a
particular affine connection, so they are consequences of \emph{purely
geometric properties}, at last, depending upon the curvature determined by
the connection upon the fibration. A natural preorder upon the set of fibers
of the second fibration is considered. Some remarks about the applicability
to economics and finance of the theories presented in the paper and about
the possible developments are made in the directions followed in papers 
\cite{car1}, \cite{car6}, \cite{car7}, \cite{car8}, \cite{car9} of the
author.
\end{abstract}

\bigskip

\bigskip

\section{Preliminaries}

\bigskip

\bigskip

For the general theory of fibrations we follow \cite{d}. A fibration or
fiber space is a pair $\mathcal{F}=(X,\pi )$, where

i) $X$ is a non-empty set, said the underlying set of the fiber space;

ii) $\pi $ is a surjection of $X$ onto a non-empty set $B$, called the base
of the fiber space;

iii) for any point $b$ in $B$ there is a subset $U$ of $B$ containing $b$, a
set $F_{b}$ and a bijection 
\[
h:U\times F_{b}\rightarrow \pi ^{-1}(U) 
\]
such that 
\[
\pi (h(y,t))=y, 
\]
for each $y$ in $U$ and $t$ in $F_{b}$. In other terms, 
\[
\pi \circ h=\mathrm{pr}_{1}^{U\times F_{b}}, 
\]
where $\mathrm{pr}_{1}^{U\times F_{b}}$ is the first projection of the
cartesian product $U\times F_{b}$.

\bigskip

Let $k$ be a natural number (an integer greater or equal to $0$) a $C^{k}$%
-fibration or fiber space of class $C^{k}$ is a pair $\mathcal{F}=(X,\pi )$,
where

i) $X$ is a $C^{k}$-manifold, said the underlying set of the fiber space;

ii) $\pi $ is a surjection of $X$ onto a $C^{k}$-differentiable manifolds $%
(B,\mathcal{A})$;

iii) for any point $b$ in $B$ there is an open neighborhood $U$ of $b$ in $B$%
, a differentiable manifold $(F,\mathcal{A}_{F})$ and a $C^{k}$%
-diffeomorphism 
\[
h:U\times F\rightarrow \pi ^{-1}(U) 
\]
such that 
\[
\pi (h(y,t))=y, 
\]
for each $y$ in $U$ and $t$ in $F$.

\bigskip

\bigskip

\section{Fibrations on financial events plane}

\bigskip

\bigskip

In this section we introduce the basic concepts of the paper.

\begin{theorem}
The space of financial events $\Bbb{R}^{2}$ is a smooth fiber space in the
following two ways:

1) the trivial one $(\Bbb{R}^{2},\mathrm{pr}_{1})$;

2) $\mathcal{F}_{i}=(\Bbb{R}^{2},\pi _{i})$ with $i>-1$ and $\pi _{i}$ the
below surjection 
\[
\pi _{i}:\Bbb{R}^{2}\rightarrow \Bbb{R}:(t,c)\mapsto (1+i)^{-t}c.
\]
\end{theorem}

\emph{Proof.} Straightforward by definition of fiber space. $\blacksquare $

\begin{definition}
We call the fibration $(\Bbb{R}^{2},\mathrm{pr}_{1})$ \textbf{natural
fibration of the financial events plane}. We call the fibration $\mathcal{F}%
_{i}=(\Bbb{R}^{2},\pi _{i})$ \textbf{fibration induced on the financial
events plane by the compound capitalization at rate} $i>-1$.
\end{definition}

Let us examine the fibration $(\Bbb{R}^{2},\mathrm{pr}_{1})$:

\begin{itemize}
\item  the base of the fibration is the time-line $\Bbb{R}$;

\item  for each time $t$, the fiber $\left( \Bbb{R}^{2}\right) _{t}$ is the
straight-line $\mathrm{pr}_{1}^{-}(t)=\left\{ t\right\} \times \Bbb{R}$,
that is the equivalence class generated by the null event $(t,0)$ by means
of the equivalence relation ``to have the same time'';

\item  this fibration is a fibred space of fiber-type $\Bbb{R}$, since each
fiber is diffeomorphic to the standard manifold $\Bbb{R}$.
\end{itemize}

Let us examine the second fibration $\mathcal{F}_{i}=(\Bbb{R}^{2},\pi _{i})$:

\begin{itemize}
\item  the base of the fibration is the capital-line;

\item  for each element $c$ of the capital-line, the fiber $(\Bbb{R}^{2})_{c}
$ is the set-curve $\pi _{i}^{-}(c)=\mathrm{gr}(M_{c})$, graph of the
function 
\[
M_{c}:\Bbb{R\rightarrow R}:t\mapsto (1+i)^{t}c,
\]
the so called capital-evolution of the event $(0,c)$. The fiber $\pi %
_{i}^{-}(c)$ is nothing but the class of equivalence generated by the event $%
(0,c)$ by means of the equivalence relation $\sim _{i}$ induced by the
compound capitalization at rate $i$, that is the binary relation defined by 
\[
e_{0}\sim _{i}e\;\;\mathrm{iff}\;\;\pi _{i}(e_{0})=\pi _{i}(e),
\]
the equivalence class generated by an event $e$ shall be denoted also by $%
\left[ e\right] _{i}$;

\item  this fibration is a fibred space of fiber-type $\Bbb{R}$, since each
fiber is diffeomorphic to the standard manifold $\Bbb{R}$ (since each fiber
is a set-curve).
\end{itemize}

\bigskip

\textbf{Remark (the fibration induced by a capitalization factor at }$0$%
\textbf{).} If $f:\left[ 0,+\infty \right] \rightarrow \Bbb{R}$ is a
capitalization factor of class $C^{k}$, that is a positive function from the
time semi-line $\left[ 0,+\infty \right] $ into the capital line $\Bbb{R}$
of class $C^{k}$ such that $f(0)=1$, we can build up a $C^{k}$-fibration $%
(\left[ 0,+\infty \right] \times \Bbb{R},\pi _{f}\Bbb{)}$, defined by 
\[
\pi _{f}:\left[ 0,+\infty \right] \times \Bbb{R}\rightarrow \Bbb{R}%
:(t,c)\mapsto f(t)^{-1}c. 
\]
Even more generally, we can define a $C^{0}$-fibration $(\Bbb{R}^{2},\pi _{f}%
\Bbb{)}$ by 
\[
\pi _{f}:\Bbb{R}^{2}\rightarrow \Bbb{R}:(t,c)\mapsto \left\{ 
\begin{array}{ccc}
f(t)^{-1}c & if & t\geq 0 \\ 
f(-t)c & if & t<0
\end{array}
\right. , 
\]
and this fibration is at least of class $C^{1}$ if $k>0$. Indeed, setting $%
g_{>}(t)=f(t)^{-1}$ and $g_{<}(t)=f(-t)$, we have $g_{<}^{\prime
}(t)=-f^{\prime }(-t)$ and $g_{>}^{\prime }(t)=-f^{\prime }(t)f(t)^{-2}$,
from which 
\[
g_{>}^{\prime }(0)=g_{<}^{\prime }(0)=-f^{\prime }(0). 
\]

\bigskip

\bigskip

\section{Properties of the fibration induced by the compound interest}

\bigskip

\bigskip

\begin{theorem}
Let, for any real $i>-1$, $\mathcal{F}_{i}=(\Bbb{R}^{2},\pi _{i})$ be the
fibration induced by the compound capitalization at rate $i$. Then, for any
two rates $i$ and $i^{\prime }$ the two fibrations $\mathcal{F}_{i}$ and $%
\mathcal{F}_{i^{\prime }}$ are isomorph, being the bijection $g:\Bbb{R}%
^{2}\rightarrow \Bbb{R}^{2}$ defined by 
\[
g(t,c)=(t,(u^{\prime })^{t}u^{-t})=(t,(u^{\prime }/u)^{t}c),
\]
for any financial event $(t,c)$, where $u=1+i$ and $u^{\prime }=1+i^{\prime }
$, an $\Bbb{R}$-isomorphism.
\end{theorem}

\emph{Proof.} An isomorphism of a $C^{k}$-fibration $\mathcal{F}=(X,\pi )$
onto another $C^{k}$-fibration $\mathcal{F}^{\prime }=(X^{\prime },\pi
^{\prime })$ with the same base $B$ is a pair of $C^{k}$-functions $(\mathrm{%
id}_{B},g)$, with $g:X\rightarrow X^{\prime }$, such that 
\[
\pi ^{\prime }\circ g=\pi . 
\]
Put $u=1+i$, $u^{\prime }=1+i^{\prime }$ and consider the bijection $g:\Bbb{R%
}^{2}\rightarrow \Bbb{R}^{2}$ defined by 
\[
g(t,c)=(t,(u^{\prime })^{t}u^{-t}c)=(t,(u^{\prime }/u)^{t}c), 
\]
for any financial event $(t,c)$, then the pair $(\mathrm{id}_{\Bbb{R}},g)$
is an isomorphism of $\mathcal{F}_{i}$ onto $\mathcal{F}_{i^{\prime }}$.
Indeed we have 
\begin{eqnarray*}
\pi _{i^{\prime }}(g(t,c)) &=&\pi _{i^{\prime }}((t,(u^{\prime
})^{t}u^{-t}c))= \\
&=&(u^{\prime })^{-t}(u^{\prime })^{t}u^{-t}c= \\
&=&\pi _{i}(t,c),
\end{eqnarray*}
for each financial event $(t,c)$. $\blacksquare $

\bigskip

\textbf{Remark.} Another way to prove that the two above induced fibrations
are isomorph is to prove that, for every $c_{0}$ belonging to the common
base $\Bbb{R}$ there is an isomorphism $g_{c_{0}}:X_{c_{0}}\rightarrow
X_{c_{0}}^{\prime }$. Indeed, define 
\[
g_{c_{0}}(t,c)=(t,(u^{\prime }/u)^{t}c), 
\]
for every event $(t,c)$ in the fiber $X_{c_{0}}=\left[ (0,c_{0})\right] _{i}$%
. We note that if $(t,c)$ is a financial event of the fiber generated by the
event $(0,c_{0})$, it has the form $(t,c_{0}u^{t})$, applying the function $%
g_{c_{0}}$ we obtain 
\begin{eqnarray*}
g_{c_{0}}(t,c) &=&(t,(u^{\prime }/u)^{t}c)= \\
&=&(t,(u^{\prime }/u)^{t}c_{0}u^{t})= \\
&=&(t,(u^{\prime })^{t}c_{0}),
\end{eqnarray*}
that is an event of the fiber $X_{c_{0}}^{\prime }=\left[ (0,c_{0})\right]
_{i^{\prime }}$: in other words, we pull back the event $e$ along the fiber $%
X_{c_{0}}$ to the event $e_{0}=(0,c_{0})$ and then we push forward the event 
$e_{0}$ to $e^{\prime }=(t,(u^{\prime })^{t}c_{0})$ along the fiber $%
X_{c_{0}}^{\prime }$.

\begin{corollary}
Let, for any real $i>-1$, $\mathcal{F}_{i}=(\Bbb{R}^{2},\pi _{i})$ be the
fibration induced by the compound capitalization at rate $i$. Then, $%
\mathcal{F}_{i}$ is trivializable for every rate $i$.
\end{corollary}

\emph{Proof.} It derives from the circumstance that the fibration $\mathcal{F%
}_{0}$ (corresponding to the rate $0\%$) is trivializable, in fact the
projection $\pi _{0}$ acts as follows 
\[
\pi _{0}:\Bbb{R}^{2}\rightarrow \Bbb{R}:(t,c)\mapsto (1+0)^{-t}c=c, 
\]
and then the projection $\pi _{0}$ is nothing but $\mathrm{pr}_{2}$ on the
cartesian product of the time-line $T$ times the capital line $C$; now it is
clear that this fibration is isomorph to the fibration $(C\times T,\mathrm{pr%
}_{1})$. The conclusion follows from the fact that each fibration $\mathcal{F%
}_{i}$ is isomorph to the fibration $\mathcal{F}_{0}$. $\blacksquare $

\bigskip

\bigskip

\section{Sections of the fibration induced by the compound interest}

\bigskip

\bigskip

\begin{theorem}
\textbf{(the sections of the fibrations }$\mathcal{F}_{i}$\textbf{).} Let $C$
be the real line of capitals and let $E$ be the plane of financial events.
Then, a curve $s:C\rightarrow E$ defined by $s(c)=(s_{1}(c),s_{2}(c))$, for
every capital $c$, is a section of the fibration $\mathcal{F}_{i}$ if and
only if 
\[
s_{2}(c)=c(1+i)^{s_{1}(c)},
\]
for every capital $c$.
\end{theorem}

\bigskip

\emph{Proof.} The curve $s$ is a section of the fibration $\mathcal{F}_{i}$,
by definition, if and only if 
\[
\pi _{i}(s(c))=c, 
\]
for every capital $c$. This last relation means that 
\[
\pi _{i}(s_{1}(c),s_{2}(c))=(1+i)^{-s_{1}(c)}s_{2}(c)=c, 
\]
for any capital $c$, that is 
\[
s_{2}(c)=c(1+i)^{s_{1}(c)}, 
\]
for any capital $c$. $\blacksquare $

\bigskip

\textbf{Remark.} In other words, the above theorem states that are sections
of the fibration induced by the compound capitalization at rate $i$ only
those curves $s:C\rightarrow E$ of the form 
\[
s(c)=\left( s_{1}(c),c(1+i)^{s_{1}(c)}\right) , 
\]
for every $c$ in $C$ and for any function $s_{1}:C\rightarrow E$.

\bigskip

We can restate the above theorem as follows.

\begin{theorem}
\textbf{(the sections of the fibrations }$\mathcal{F}_{i}$\textbf{).} Let $C$
be the capital line, $T$ the time line and let $E$ be the plane of financial
events. Then, a curve $s:C\rightarrow E$ is a section of the fibration $%
\mathcal{F}_{i}=(E,\pi _{i})$ if and only if there if a function $%
f:C\rightarrow T$ such that 
\[
s(c)=(f(c),u^{f(c)}c),
\]
for every capital $c$.
\end{theorem}

\bigskip

\textbf{Remark.} The fibration $\mathcal{F}_{0}$ is the pair $(E,\mathrm{pr}%
_{2})$, thus every its section $s:C\rightarrow E$ has the form 
\[
s(c)=(f(c),c), 
\]
for every $c$ in $C$, where $f:C\rightarrow T$ is any function of the
capital line into the time line. Since any fibration $\mathcal{F}_{i}$ is
isomorphic to the fibration $\mathcal{F}_{0}$, the sections of $\mathcal{F}%
_{i}$ can be obtained by the section of $\mathcal{F}_{0}$ applying the
canonical isomorphism of $\mathcal{F}_{0}$ into $\mathcal{F}_{i}$, that is
the bijection $g:\Bbb{R}^{2}\rightarrow \Bbb{R}^{2}$ defined by 
\[
g(t,c)=(t,u^{t}c)=(t,u^{t}c), 
\]
for any financial event $(t,c)$, where $u=1+i$; applying the isomorphism $g$
to the section $s$, we obtain the curve $g\circ s:C\rightarrow E$, that is
the curve defined by 
\[
g\circ s(c)=g(f(c),c)=(f(c),u^{f(c)}c), 
\]
for any capital $c$: so we obtained newly the above theorem.

\bigskip

In a perfectly analogous way we can extend the above theorem as follows.

\begin{theorem}
\textbf{(the sections of the fibrations }$\mathcal{F}_{i}$\textbf{).} Let $%
C^{\prime }$ be a part of the capital line, $T^{\prime }$ be a part of the
time line and let $E$ be the plane of financial events. Then, a curve $%
s:C^{\prime }\rightarrow E$ is a section of the fibration $\mathcal{F}%
_{i}=(E,\pi _{i})$ upon the part $C^{\prime }$ if and only if there if a
function $f:C^{\prime }\rightarrow T^{\prime }$ such that 
\[
s(c)=(f(c),u^{f(c)}c),
\]
for every capital $c$ in $C^{\prime }$.
\end{theorem}

\bigskip

\bigskip

\section{Capital evolutions as sections in the compound interest}

\bigskip

\bigskip

We devote this paragraph to solve this problem important in the applications:

\begin{itemize}
\item  Let $M:T\rightarrow C$ be a function from the time line into the
capital line, called a capital evolution. There are sufficient conditions to
assure that the graph of the function $M$, the subset $\mathrm{gr}(M)$ of
the financial events plane $E$, is the trace of a section $s:C\rightarrow E$?
\end{itemize}

\bigskip

At this purpose we have the following complete result.

\begin{theorem}
Let $M:T\rightarrow C$ be a function from the time line into the capital
line. Then, the graph of the function $M$, the subset $\mathrm{gr}(M)$ of $E$%
, is the trace of a section $s:C\rightarrow E$ if and only if there exists a
bijection $f:C\rightarrow T$ such that 
\[
M(t)=u^{t}f^{-}(t),
\]
for each time $t$.
\end{theorem}

\emph{Proof.} \textbf{Sufficiency.} Let us suppose there exist a bijection $%
f:C\rightarrow T$ such that 
\[
M(t)=u^{t}f^{-}(t), 
\]
for each time $t$. Then, let $s:C\rightarrow E$ be the curve defined by 
\[
c\mapsto (f(c),u^{f(c)}c), 
\]
for every capital $c$. For each $t$ in $T$ (by surjectivity of the function $%
f$) there is a capital $c$ in $C$ such that $f(c)=t$, hence we have 
\[
(t,M(t))=(f(c),u^{t}f^{-}(t))=(f(c),u^{f(c)}c)=s(c), 
\]
so any point $(t,M(t))$ of the graph of $M$ is a point of the curve $s$,
that is 
\[
\mathrm{gr}(M)\subseteq s(C); 
\]
we have now to prove that $s(C)\subseteq \mathrm{gr}(M)$, indeed, let $c$ be
a capital, then by surjectivity of the reciprocal function $f^{-}$, there is
a time $t$ such that $f^{-}(t)=c$, now 
\[
s(c)=(f(c),u^{f(c)}c)=(f(c),u^{t}f^{-}(t))=(t,M(t)), 
\]
as we desire. \textbf{Necessity.} Suppose now that the graph of $M$ is the
trace of a section $s$, this is equivalent to say (by the above
characterization of sections) that there is a function $f:C\rightarrow T$
(not necessarily a bijection) such that, for each time $t$ in $T$, we have 
\[
(t,M(t))=(f(c),u^{f(c)}c), 
\]
for some $c$ in $C$. First of all, we have to prove that the function $f$ is
bijective. In fact, let $c$ and $c^{\prime }$ be two capitals such that $%
f(c)=f(c^{\prime })$, since $f(c)$ is in $T$, we have 
\begin{eqnarray*}
(f(c),M(f(c)) &=&s(c)=(f(c),u^{f(c)}c), \\
(f(c^{\prime }),M(f(c^{\prime })) &=&s(c^{\prime })=(f(c^{\prime
}),u^{f(c^{\prime })}c^{\prime }),
\end{eqnarray*}
from which 
\begin{eqnarray*}
u^{f(c)}c &=&M(f(c))= \\
&=&M(f(c^{\prime }))= \\
&=&u^{f(c^{\prime })}c^{\prime }= \\
&=&u^{f(c)}c^{\prime },
\end{eqnarray*}
and we conclude $c=c^{\prime }$. The function $f$ is then injective, it is
surjective since for every $t$ there is a $c$ such that $t=f(c)$. Concluding
the relation $M(t)=u^{t}f^{-}(t)$, is an obvious consequence of the
relations $t=f(c)$ and $M(t)=u^{f(c)}c$ by means of bijectivity. The theorem
is proved. $\blacksquare $

\bigskip

We conclude the section with a little (sometimes useful) result.

\bigskip

\begin{proposition}
Let $M:T\rightarrow C$ be a capital evolution. Then, the graph of $M$ is the
trace of a section $s:C\rightarrow E$ of the fibration $\mathcal{F}_{i}$ if
and only if the mapping 
\[
h:\mathrm{gr}(M)\rightarrow C:(t,c)\mapsto cu^{-t}
\]
is a bijection.
\end{proposition}

\emph{Proof.} \textbf{Necessity.} Let the graph of $M$ be a section, then
there is a bijection $f:C\rightarrow T$ such that 
\[
M(t)=u^{t}f^{-}(t), 
\]
for each time $t$. We have 
\begin{eqnarray*}
h(t,c) &=&h(t,M(t))= \\
&=&h(t,u^{t}f^{-}(t))= \\
&=&u^{t}f^{-}(t)u^{-t}= \\
&=&f^{-}(t),
\end{eqnarray*}
and this prove that $h$ is a bijection. \textbf{Sufficiency.} Let the
mapping $h$ be a bijection, we put 
\[
v(t)=h(t,M(t)), 
\]
it is clear that $v$ is a bijection, moreover 
\[
v(t)=h(t,M(t))=M(t)u^{-t}, 
\]
from which, setting $f=v^{-}$, we deduce, for each $t$ in $T$, 
\[
M(t)=v(t)u^{t}=f^{-}(t)u^{t}, 
\]
as we desired. $\blacksquare $

\bigskip

Analogous result we have for the evolutions defined on a part of the time
line.

\begin{theorem}
Let $T^{\prime }$ be a part of the time line, $C^{\prime }$ be a part of the
capital line and $M:T^{\prime }\rightarrow C^{\prime }$ be a capital
evolution. Then, the graph of the function $M$, that is a subset $\mathrm{gr}%
(M)$ of the rectangle $T^{\prime }\times C^{\prime }$, is the trace of a
section $s:C^{\prime }\rightarrow E$ upon the part if and only if there
exists a bijection $f:C^{\prime }\rightarrow T^{\prime }$ such that 
\[
M(t)=u^{t}f^{-}(t),
\]
for each time $t$ in $T^{\prime }$.
\end{theorem}

\bigskip

\textbf{Example.} Let $i$ be a positive real, $T$ be the time line, let $%
C_{>}$ be the semi-line of strictly positive capital and let $M:T\rightarrow
C_{>}$ be a surjective $C^{1}$-capital evolution such that $M^{\prime }$ is
strictly negative. Then, the graph of $M$ is the trace of a section of the
fibration $\mathcal{F}_{i}$ upon $C_{>}$. Indeed, put 
\[
v(t)=M(t)u^{-t}, 
\]
we have 
\[
v^{\prime }(t)=M^{\prime }(t)u^{-t}-M(t)u^{-t}\ln u<0, 
\]
for any time $t$, so the function $v$ is strictly decreasing (hence
injective) and surjective since $M$ is so, and the claim is proved taking
for $f$ the inverse of $v$.

\bigskip

\bigskip

\section{Connections on the financial fibration and capitalization laws}

\bigskip

\bigskip

Consider the trivial financial fibration $(E,\mathrm{pr}_{1})$, where $E$ is
the rectangle $U\times \Bbb{R}$ product of an open subset of the time-line
times the capital axis $\Bbb{R}$.

\begin{definition}
\textbf{(of local discount factor). }Let $t$ be a time in $U$ and let $%
F:V\rightarrow \Bbb{R}$ be an application of a neighborhood $V$ of the
time-vector $0$ into the discount factor line $\Bbb{R}$. The mapping $F$ is
said a (local) discount law in $U$ if verifies the following properties:

i) the translation $t+V$ is contained in the open subset $U$;

ii) the discount factor $F(h)$ is positive, for every time-vector $h$;

iii) the $F$-discount factor at time $0$ is $1$, $F(0)=1$;

i) $F$ is of class $C^{1}$ in $V$.
\end{definition}

\begin{definition}
\textbf{(of financial translation induced by a discount factor). }We call,
for every time-vector $h$ in $\Bbb{R}$, such that $t+h$ lies in $U$, \textbf{%
financial translation from the fiber} $E_{t}$ \textbf{to the fiber} $E_{t+h}$
\textbf{induced by the discount law} $F$ the mapping 
\[
\tau _{h}:\left\{ t\right\} \times \Bbb{R}\rightarrow \left\{ t+h\right\} 
\times \Bbb{R}
\]
defined by 
\[
\tau _{h}:(t,c)\mapsto (t+h,F(h)^{-1}c),
\]
for every financial event $e=(t,c)$ of the fiber $E_{t}$.
\end{definition}

\begin{theorem}
The financial translation $\tau _{h}$ induced by a discount law $F$ is a
linear isomorphism of the fiber $E_{t}=\left\{ t\right\} \times \Bbb{R}$
onto the fiber $E_{t+h}=\left\{ t+h\right\} \times \Bbb{R}$ and the
application $\tau $ of $V\times \Bbb{R}$ into $U\times \Bbb{R}$ defined by 
\[
\tau :(h,c)\mapsto (t+h,F(h)^{-1}c)
\]
is of class $C^{1}$. The derivative $\tau ^{\prime }(0,c)$ of the
application $\tau $ at the point $(0,c)$ is the linear mapping of $\Bbb{R}%
\times \Bbb{R}$ into itself 
\[
(k,v)\mapsto (k,v-F^{\prime }(0)kc).
\]
\end{theorem}

\begin{theorem}
Let $F$ be a discount law. Then, the mapping $(k,c)\mapsto F^{\prime }(0)kc$
is a bilinear application of $\Bbb{R}\times \Bbb{R}$ into $\Bbb{R}$, we
denote it by $\Gamma _{t}$ (and we call it the Cristoffel bilinear form) 
\[
(k,c)\mapsto \Gamma _{t}(k,c)=F^{\prime }(0)kc.
\]
Conversely, if we have a bilinear application $(k,c)\mapsto \Gamma _{t}(k,c)$
and if we put 
\[
F(h)=1+\Gamma _{t}(h,1),
\]
the function $F$ is a discount factor such that 
\[
F^{\prime }(0)kc=\Gamma _{t}(k,c).
\]
\end{theorem}

Since $E=U\times \Bbb{R}$ and since the event $e=(t,c)$ is a point of a
fiber $E_{t}$, the tangent space $T_{(t,c)}(E)$ can be identified with the
product $T_{t}(U)\times T_{c}(\Bbb{R})$, and this product can be itself
identified with the product $(\left\{ t\right\} \times \Bbb{R})\times %
(\left\{ c\right\} \times \Bbb{R})$.

\begin{definition}
\textbf{(of local connection induced by a discount factor).} Let $\mathcal{T}
$ be the time line endowed with its natural structure of $C^{\infty }$
manifold. We call the application $C_{t}$ of the product $T_{t}(\mathcal{T})%
\times E_{t}$ into the tangent bundle $T(E)$ of the fibration $(E,\mathrm{pr}%
_{1})$, union of the (disjoint) tangent spaces $T_{e}(E)=\left\{ e\right\} 
\times \Bbb{R}^{2}$, with $e$ varing in $E$, defined by 
\[
C_{t}:T_{t}(\mathcal{T})\times E_{t}\rightarrow T(E):((t,k),e)\mapsto
(e,(k,-F^{\prime }(0)kc)),
\]
local connection at time $t$ induced by the discount law $F$. The local
connection $C_{t}$ associate with a couple of (applied) vectors $(t,k)\in
T_{t}(\mathcal{T})$, $(t,c)\in E_{t}$ the applied vector at the event $%
e=(t,c)$ given by
\[
C_{t}((t,k),(t,c))=((t,c),(k,-\Gamma _{t}(k,c))).
\]
\end{definition}

\begin{definition}
\textbf{(of global connection induced by a discount factor). }Let $F$ be a
global discount law. The connection induced by the capitalization factor $F$
is the mapping 
\[
C:T(\mathcal{T})\oplus E\rightarrow T(E):((t,k),e)\mapsto C_{t}((t,k),e),
\]
where $T(\mathcal{T})\oplus E$ is the union of the (disjoint) rectangles $%
T_{t}(\mathcal{T})\times E_{t}$, i.e., the rectangles $\left( \left\{
t\right\} \times \Bbb{R}\right) \times \left( \left\{ t\right\} \times \Bbb{R%
}\right) $.
\end{definition}

\textbf{Application.} Consider an event $e=(t,c)$ and a capitalization law $%
u:\Bbb{R}\rightarrow \Bbb{R}$, that is a mapping verifing the following
properties:

\begin{itemize}
\item  the capitalization factor $u(h)$ is positive, for every time-vector $h
$;

\item  the $u$-capitalization factor at time $0$ is $1$, i.e. $u(0)=1$;

\item  $u$ is of class $C^{1}$
\end{itemize}

The capital-evolution of the event $e$ determined by the capitalization
factor $u$ is by definition the mapping 
\[
M:T\rightarrow C:t\mapsto u(t)c.
\]
We note that the moltiplicative inverse $v=u^{-1}$ is a discount law. Let us
consider the connection induced by the discount factor $v$: 
\[
C_{t}:T_{t}(\mathcal{T})\times E_{t}\rightarrow T(E):((t,k),e)\mapsto
(e,(k,-F^{\prime }(0)kc)).
\]
Suppose that each event $e=(t,c)$ has a capitalization-time $t$, that is we
suppose that $e$ is the state at $t$ of the event $e_{0}=(0,cu(t)^{-1})$,
the financial translation induced by the capitalization law $u$ is defined
by 
\[
\tau _{h}:(t,c)\mapsto (t+h,u(t+h)u(t)^{-1}c),
\]
so, concerning the discount law we have
\[
v(h)^{-1}=u(t+h)u(t)^{-1},
\]
deriving we obtain
\[
-v(h)^{-2}v^{\prime }(h)=u^{\prime }(t+h)u(t)^{-1},
\]
and considering the Cristoffel bilinear form, we have
\begin{eqnarray*}
-\Gamma _{t}(k,c) &=&-v^{\prime }(0)kc= \\
&=&u^{\prime }(t)u(t)^{-1}kc= \\
&=&\delta (t)kc,
\end{eqnarray*}
where $\delta (t):=u^{\prime }(t)u(t)^{-1}$ is the instant force of interest
(by definition) at time $t$ of the capitalization law $u$.

 \end{document}